\documentclass[11pt]{article}

\usepackage{graphicx}

\usepackage{amssymb,subfigure}
\usepackage{amsmath}

\newtheorem{teo}{Theorem}

\newtheorem{rem}{Remark}
\newtheorem{lem}{Lemma}
\newtheorem{prop}{Proposition}

\renewcommand{\r}{{\mathbb R}}
\newcommand{\z} {{\mathbb Z}}
\newcommand{\cn} {{\mathbb C}}
\newcommand{\n} {{\mathbb N}}

\newcommand{\be}{\begin{equation}}
\newcommand{\ee}{\end{equation}}
\newcommand{\bega}{\begin{gather}}
\newcommand{\enga}{\end{gather}}

\begin{document}

\title{On a connection between nonstationary and periodic wavelets
}

\author{Elena A. Lebedeva\footnote{Mathematics and Mechanics Faculty, Saint Petersburg State University,
Universitetsky prospekt, 28, Peterhof,  Saint Petersburg,
 198504, Russia; Saint Petersburg State Polytechnical University,
 Polytechnicheskay 29,  Saint Petersburg, 195251, Russia}}

\date{ealebedeva2004@gmail.com}

\maketitle


\begin{abstract}
We compare  frameworks of nonstationary nonperiodic wavelets and periodic wavelets. We construct one system from another using periodization.   There are infinitely many nonstationary systems corresponding to the same periodic wavelet. Under mild conditions on periodic scaling functions, among these nonstationary wavelet systems, we find a system such that its time-frequency localization is adjusted with an angular-frequency localization of an initial periodic wavelet system. Namely, we get the following equality 
$
\lim_{j\to \infty} UC_B(\psi^P_j) = \lim_{j\to \infty}UC_H(\psi^N_j),
$           
where $UC_B$ and $UC_H$ are the Breitenberger and the Heisenberg uncertainty constants, $\psi^P_j \in L_2(\mathbb{T})$ and 
$\psi^N_j\in L_2(\mathbb{R})$  are periodic and nonstationary wavelet functions respectively.   
\end{abstract}



\section{Introduction}
\label{intr}

The frameworks of nonstationary nonperiodic and periodic wavelets are introduced and studied separately. 
The notion of nonstationary wavelet system is introduced independently by M.~Z.~Berkolayko, I.~Y.~Novikov  \cite{BN} and by C.~de Boor, R.~DeVore, A.~Ron \cite{BDR}. In \cite{BN}, the nonstationary system (called almost-wavelets) is used to construct  
an orthonormal shift invariant basis consisting of infinitely differentiable compactly supported functions. It is well known that it is impossible to construct stationary wavelet basis satisfying these properties. Further, nonstationary wavelets are studied in 
\cite{DR, DKLR, HSh, NPS}. 
Concerning to periodic case, first, periodic wavelets are generated by periodization of stationary wavelet functions. Later, a general approach to study periodic wavelets directly using a periodic analog of multiresolution analysis (MRA) is developed. A notion of periodic MRA is introduced and discussed in 
\cite{klt95,  NPS, P, PT, S97}.

The main purpose of this paper is to compare   frameworks of nonstationary nonperiodic wavelets and periodic wavelets. We intend to extrapolate the idea of periodization and to study how to construct one system from another using periodization. Parseval wavelet frames generated by unitary extension principles (UEP) are considered. We have not found a UEP for nonstationary setting in literature, so we formulate this principle in Theorem \ref{nUEP}. The proof is omitted since it can be easily rewritten from  Theorem 1.8.9 \cite{NPS} replacing stationary conditions by nonstationary ones. For periodic setting we cite UEP proved by  S.S.~Goh,  K.M.~Teo in 2008 (Theorem \ref{pUEP}). 

In Proposition \ref{nst_to_per} it  is proved that the periodization of a nonstationary Parseval wavelet frame is a periodic Parseval wavelet frame. The proof is reduced to comparison of the UEPs. The inverse problem, namely to construct a nonstationary system such that its periodization is the initial periodic system, is more interesting and complicated. It has infinitely many solutions. In Lemma \ref{def_mask} we design a family of nonstationary masks $m^K_j$  corresponding to the initial periodic wavelets. A family parameter $K\in\n$ is responsible for the smoothness of a mask and for an order of a zero at the point $\xi=1/2.$ It is well known that the last characteristic is also important since it is necessary condition for the smoothness of scaling and wavelet function, at least, for stationary case. In Lemma \ref{prodL2}, we suggest a nonstationary analogue of a sufficient condition for an infinite product 
$\prod_{r=1}^{\infty} m_{j+r}(\xi 2^{-r})$ to converge uniformly on a compact and to be in $L_2(\r).$ 
In Lemma \ref{my_prodL2}, this sufficient condition is specified for the masks $m^K_j$ defined in Lemma \ref{def_mask}. So, we can design nonstationary wavelet system such that its periodization coincides with the given periodic system and its nonstationary masks can be chosen arbitrarily smooth. 

Now, among these solutions we are looking for wavelet frame with the following additional property. Time-frequency localization of the resulting nonstationary system should be adjusted with angular-frequency localization of the initial periodic system, that is 
$
\lim_{j\to \infty} UC_B(\psi^P_j) = \lim_{j\to \infty}UC_H(\psi^N_j),
$  
where $UC_B,$ $UC_H$ are the Breitenberger and the Heisenberg uncertainty constants (UCs).      
 We have a particular motivation to be interested in this issue. Good time-frequency localization (that means finiteness of UC) is a desirable and natural  property of wavelet functions. It follows from the uncertainty principle that there is no function with arbitrarily small UC (less than $1/2$ for both UCs). In \cite{Bat} it is proved that there is no wavelet function with $UC_H$ less than $3/2.$  And it is unknown whether there exists a stationary or a nonstationary wavelet (an orthonormal basis or a tight frame) $\psi^0$ such that 
$UC_H(\psi^0) \leq  3/2 +\varepsilon.$ 
Results on  (un)boundedness of UC for various wavelet families can be found  in
\cite{GB, GL, L07, L12, L11, N1}. However in the periodic case the situation is little bit better. In \cite{LebPres14} we construct  a Parseval periodic wavelet frame satisfying the following property
$\lim_{j\to\infty} UC_B(\psi^{P}_j)=3/2$. And in \cite{L15}, we show that $3/2$ is a minimal possible value of $UC_B$ for a wide class of periodic wavelet sequences. 
So, our particular motivation is to construct  nonstationary wavelet systems with good localization  starting with well-localized periodic ones, and conversely.   
In Theorem \ref{nstPQRS} we obtain a sufficient condition for nonstationary and periodic wavelet  frames to satisfy equality 
$
\lim_{j\to \infty} UC_B(\psi^P_j) = \lim_{j\to \infty}UC_H(\psi^N_j).
$   
In Theorem \ref{adjust_loc} this sufficient condition is formulated in terms of scaling masks of the given periodic wavelet system.

\section{Notations and auxiliary results}
\label{note}
Let $L_2(\mathbb{T})$ be the space of all $1$-periodic square-integrable complex-valued functions, with
inner product $(\cdot,\cdot)$ given by 
$
(f,\,g):=  \int_{0}^{1} f(x)\overline{g(x)}\,\mathrm{d}x
$
for any $f,g \in L_2(\mathbb{T}),$ and norm $\|\cdot\|:=\sqrt{(\cdot,\,\cdot)}.$ 
The Fourier series of a function 
$
f \in L_2(\mathbb{T})
$
is defined by
$\sum_{k \in \mathbb{Z}}\widehat{f}(k) \mathrm{e}^{2 \pi \mathrm{i} k x},$
where its Fourier coefficient is defined by
$
\widehat{f}(k) = \int_{0}^{1} f(x)\mathrm{e}^{- 2 \pi \mathrm{i} k x}\,\mathrm{d}x.
$

Let $L_{2}(\r)$ be the space of all square-integrable complex-valued functions, with
inner product $(\cdot,\cdot)$ given by 
$
(f,\,g):=  \int_{\r} f(x)\overline{g(x)}\,\mathrm{d}x
$
for any $f,g \in L_2(\r),$ and norm $\|\cdot\|:=\sqrt{(\cdot,\,\cdot)}.$ 
The Fourier transform of a function 
$
f \in L_{2}(\r)
$
is defined by
$
\widehat{f}(\xi):=  \int_{\r}f(x) \mathrm{e}^{- 2 \pi \mathrm{i} \xi x}\,\mathrm{d}x.
$


We recall the notion of a tight frame.
Let $H$ be a separable Hilbert space. If there exists a constant $A>0$
such that for any $f \in H$ the following equality holds
$
\sum_{n=1}^{\infty} \left|(f,\,f_n)\right|^2 = A \|f\|^2,
$ 
 then the sequence $(f_n)_{n \in \mathbb{N}}$ is called \texttt{a tight frame} for $H.$ In the case $A=1$, a tight frame is called \texttt{ a Parseval frame}.
In addition, if $\|f_n\|=1$ for all $n \in \mathbb{N}$,  then a Parseval frame  forms an orthonormal basis.


To keep ideas clear and notations simple we consider one-dimensional wavelet systems with one wavelet generator. We recall the basic definitions.   
In the sequel,  we use the following notation for a shift  
$
f_{j,k}(x):=f_j(x-  2^{-j} k)
$
of a function $f_j \in L_2(\mathbb{T})$ or $f_j \in L_2(\mathbb{R}),$ where $k \in \mathbb{Z}.$ 
Given $\varphi^N_0,\,\psi^N_j \in L_2(\mathbb{R}),$ $j\in\z_+$, if the set 
$\Psi^N:=\left\{\varphi^N_0,  \psi^N_{j,k} : \ j\in\z_+,\ k \in \mathbb{Z} \right\}$
forms a frame (or a basis) for $L_2(\mathbb{R})$ then $\Psi^N$ is called \texttt{a nonstationary wavelet frame 
(or a nonstationary wavelet basis)} for $L_2(\mathbb{R}).$
In the stationary case, all the wavelet functions $\psi^N_j$ are generated by a single function $\psi$, namely  $\psi^N_j(x)=2^{j/2}\psi(2^j x)$.  

Analogously, let   $\varphi^P_0,\,\psi^P_j \in L_2(\mathbb{T}),$ $j\in\z_+$ be periodic functions.
 If the set 
 $\Psi^P:=\left\{\varphi^P_0,  \psi^P_{j,k} : \ j\in\z_+,\ \right.$
$\left. k=0,\dots,2^j-1 \right\}$
forms a tight frame (or a basis) for $L_2(\mathbb{T})$ then $\Psi^P$ is said to be \texttt{a periodic tight wavelet frame 
(or a periodic wavelet basis)} for $L_2(\mathbb{T}).$


We compare framework of nonstationary and periodic wavelet systems generated by unitary extension principles (UEP). It is not difficult to revise a proof of stationary UEP (see, for example, Theorem 1.8.9 \cite{NPS}) and to obtain its nonstationary counterpart. The result is formulated in the following theorem. We omit  the proof since it can be rewritten from Theorem 1.8.9 \cite{NPS} replacing stationary conditions by nonstationary ones.   

\begin{teo}[the unitary extension principle for a nonstationary setting]
\label{nUEP}
Let $\varphi^N_j \in L_{2}(\mathbb{R}),$ $j\in\z_+$
be  a sequence of functions such that 
\begin{equation}
\label{ncon1} 
    \lim_{j \to \infty}2^{j/2} \widehat{\varphi^N_j}(\xi) = 1, \ \ \ \ \ \ \  \xi \in \r.
\end{equation}

Let $m^N_{0,j} \in L_{2}(\mathbb{T})$ be a $1$-periodic function such that 
\begin{equation}
\label{ncon2}
\left|m^N_{0,j}(\xi)\right|^2 +  \left|m^N_{0,j}(\xi+1/2)\right|^2 =2,
\end{equation}
and 
\begin{equation}
\label{ncon3}
	\widehat{\varphi^N_j}(\xi)=m^N_{0,j+1}(\xi/2^{-j-1}) \widehat{\varphi^N_{j+1}}(\xi).
\end{equation}

Let $\psi^N_j \in L_{2}(\mathbb{R})$, $j\in\z_+$ 
be a sequence of functions defined by
\begin{equation}
\label{ncon4}
		\widehat{\psi^N_j}(\xi):=m^N_{1,j+1}(\xi/2^{-j-1}) 
		\widehat{\varphi^N_{j+1}}(\xi),
	\end{equation}
where $m^N_{1,j}(\xi)= e^{2 \pi i \xi} \overline{m^N_{0,j}(\xi+1/2)}. $   

Then the family
$\Psi^N  =\{\varphi^{N}_{0,k}, \, \psi^{N}_{j,k} : \  j\in \mathbb{Z}_+,\  k \in \mathbb{Z}\} $
forms a Parseval wavelet frame for $L_{2}(\mathbb{R}).$
\end{teo}
The functions $\varphi^N_j,$ $\psi^N_j,$ $m^N_{0,j},$ and $m^N_{1,j}$ are 
called \texttt{a scaling function, a wavelet function, a scaling mask and a wave\-let mask}   respectively.

\begin{teo}[\cite{GT1}; the unitary extension principle for a periodic setting]
\label{pUEP}
Let $\varphi^P_j \in L_{2}(\mathbb{T}),$ $j\in\z_+$,
be  a sequence of $1$-periodic functions such that 
\begin{equation}\label{con1}
     \lim_{j \to \infty}2^{j/2} \widehat{\varphi^P_j}(k) = 1 , \ \ \ \ \ \ \  k \in \z.
\end{equation}
Let $\mu^j_k \in \cn,$ $j\in\z_+,$ $k\in\z$, be a two-parameter sequence    such that $\mu^j_{k+2^j}=\mu^j_{k},$ 
\begin{equation}
\label{con2}
\left|\mu^j_{k}\right|^2 +  \left|\mu^j_{k+2^{j-1}}\right|^2 =2,
 \end{equation}
and
\begin{equation}
\label{con3}
	\widehat{\varphi^P_j}(k)=\mu^{j+1}_{k} 
		\widehat{\varphi^P_{j+1}}(k). 
\end{equation}
Let $\psi^P_j \in L_{2}(\mathbb{T})$, $j\in\z_+$ 
be a sequence of functions defined by
\begin{equation}
\label{con4}
		\widehat{\psi^P_j}(k):=\lambda^j_{k} 
		\widehat{\varphi^P_{j+1}}(k),
	\end{equation}
where $\lambda^j_{k}= e^{2 \pi i k/2^{-j}} \overline{\mu^j_{k+2^{j-1}}}. $   
Then the family
$\Psi^P  =\{\varphi^{P}_{0}, \, \psi^{P}_{j,k} : \  j\in \mathbb{Z}_+,
\ 		  k=0,\dots,2^j-1 \}$
forms a Parseval wavelet frame for $L_{2}(\mathbb{T}).$ 
\end{teo}
The functions $\varphi^P_j,$ $\psi^P_j$ and sequences $(\mu^j_k)_{k\in\z},$  $(\lambda^j_k)_{k\in\z}$ are 
called \texttt{a scaling function, a wavelet function, a scaling mask, and a wave\-let mask}   respectively.

It is convenient to construct both systems starting with scaling masks. 
Namely, let $\nu^{j}_{k}$ be  a sequence such that  $\nu^{j}_{k}=\nu^{j}_{k+2^j}$.	We define 
$\widehat{\xi}_{j}(k):=\prod_{r=j+1}^{\infty}\nu^{r}_{k}.$ 
If the above infinite product converges, then the scaling function, scaling mask, wavelet mask, and wavelet function in a periodic setting are defined 
	 respectively as
\begin{gather}
        \widehat{\varphi_{j}}(k):=2^{-j/2}\widehat{\xi}_{j}(k),\qquad\quad 
	\mu^{j}_k:=\sqrt{2} \nu^{j}_k,\ 
	\notag\\
	\lambda^{j}_k:=e^{2\pi i 2^{-j}k}\mu^{j}_{k+2^{j-1}},\qquad \quad
	{\widehat{\psi}}_{j}(k):=\lambda^{j+1}_{k} \widehat{\varphi}_{j+1}(k).
\notag
\end{gather}

Analogously,
suppose $m_j$ is a $1$-periodic function, 
$\widehat{\varphi_j}(\xi):=\prod_{r=1}^{\infty}m_{j+r}(\xi/2^r)$ is well-defined, then the scaling function, scaling mask, and wavelet function in a nonstationary setting are defined 
	 respectively as
\begin{gather}
\label{aux_n}
        \widehat{\varphi^N_j}(\xi) = 2^{-j/2}\widehat{\varphi_j}(\xi/2^j),\qquad\quad 
	m_j(\xi)=2^{-1/2} m^N_{0,j}(\xi),\ 
	\\
	\widehat{\psi^N_j} (\xi) = 2^{-j/2}  m_{1,j+1}(\xi/2^{j+1}) \widehat{\varphi_{j+1}}(\xi/2^{j+1})=2^{-j/2} \widehat{\psi_j}(\xi/2^j)
\notag
\end{gather}

Auxiliary scaling and wavelet functions $\varphi_j$ and $\psi_j$ are connected with nonstationary  
scaling and wavelet functions as
$\psi^N_j(x)=2^{j/2}\psi_j(2^j x)$ and
$\varphi^N_j(x)=2^{j/2}\varphi_j(2^j x).$ In stationary case $\varphi_j$ and $\psi_j$ are coincided with scaling and wavelet functions 
$\varphi$ and $\psi$.


Comparing these two UEPs we get the following
\begin{prop}
\label{nst_to_per}
If $\Psi^N = \{\varphi^{N}_{0,k}, \, \psi^{N}_{j,k} \}_{j\in \mathbb{Z}_+,\, k \in \mathbb{Z}}$ is a Parseval wavelet frame generated by UEP, 
$\varphi^{N}_{0}, \psi^{N}_{j} \in L_1(\mathbb{R}),$
 and 
$$
\psi^P_j(x):=\sum_{n\in\mathbb{Z}} \psi^N_j (x-n),
$$
then $\Psi^P =\{\varphi^{P}_{0}, \, \psi^{P}_{j,k}  \}_{j\in \mathbb{Z}_+,\, k=0,\dots,2^j-1}$ is a Parseval wavelet frame.
\end{prop}

\textbf{Proof.}
Since $\varphi^{N}_{0}, \psi^{N}_{j} \in L_1(\mathbb{R}),$ periodization
$
\psi^P_j(x)=\sum_{n\in\mathbb{Z}} \psi^N_j (x-n)
$
is well-defined. In the Fourier domain the last equality is rewritten as 
$\widehat{\psi^N_j}(k)=\widehat{\psi^P_j}(k).$ 
Therefore, substituting $k\in \z$ or $k 2^{-j},$ where $k \in \z,$ $j\in\z_+$ 
instead of $\xi\in \r$ in conditions (\ref{ncon1})-(\ref{ncon4}) we immediately get
(\ref{con1})-(\ref{con4}).  
Proposition \ref{nst_to_per} is proved. \hfill $\Diamond$

 
We recall the definitions of the  UCs and the uncertainty principles.
\texttt{The Heisenberg UC} of $f \in L_2(\mathbb{R})$ is the functional 
$UC_H(f):=\Delta(f)\Delta(\widehat{f})$ such that
$$
\Delta^2(f):=
\frac{1}{\|f\|^{2}} \int_{\mathbb{R}} \bigl(x-c(f)\bigr)^2 \bigl|f(x)\bigr|^2 \, dx , \quad   
c(f):=
\frac{1}{\|f\|^{2}} \int_{\mathbb{R}} x \bigl|f(x)\bigr|^2 \, dx,
$$
where 
$\Delta(f),$ $\Delta(\widehat{f}),$ $c(f),$ and  $c(\widehat{f})$
are called \texttt{time variance, frequency variance, time centre,} and \texttt{frequency centre} respectively. 
The Heisenberg uncertainty principle says that $UC_H(f)\geq 1/2$
for $f \in L_2(\mathbb{R})$,   and the equality is attained iff $f$ is the Gaussian function. 
In \cite[p. 137]{Bat} the following  refinement of the Heisenberg uncertainty principle is proved.
If $f \in L_2(\r),$ $\ c(\widehat{f})$ $=0$, and $\int_{\r} f = 0,$  
then $UC_H(f)\geq 3/2.$

The UC for periodic functions is introduced in \cite{B}. 
Let $f =\sum_{k \in \mathbb{Z}} c_k \mathrm{e}^{ \mathrm{i} k \cdot}\in L_{2,2\pi}.$  
\texttt{ The first trigonometric moment} is defined as
$$ 
\tau(f):=\frac{1}{2 \pi} \int_{-\pi}^{\pi} \mathrm{e}^{ \mathrm{i} x} |f(x)|^2\, \mathrm{d}x =
 \sum_{k \in \mathbb{Z}} c_{k-1} \overline{c_{k}}.
$$
\texttt{The angular variance} of the function $f$ is defined by  
$$
{\rm var_A }(f):= \frac{\left(\sum_{k \in \mathbb{Z}}|c_k|^2\right)^2}{
\left|\sum_{k \in \mathbb{Z}}c_{k-1} \bar{c}_{k}\right|^2}-1
=
\frac{\|f\|^4}{|\tau(f)|^2}-1.
$$
\texttt{The frequency variance} of the function $f$ is defined by  
$$
{\rm var_F }(f):= \frac{\sum_{k \in \mathbb{Z}}k^2 |c_k|^2}{\sum_{k \in \mathbb{Z}}|c_k|^2}-
\frac{\left(\sum_{k \in \mathbb{Z}}k|c_k|^2\right)^2}{\left(\sum_{k \in \mathbb{Z}}|c_k|^2\right)^2}
=
\frac{\|f'\|^2}{\|f\|^2}-\frac{(i f',\, f)^2}{\|f\|^4}.
$$
The quantity 
$
UC_B(\{c_k\}):=UC_B(f):=\sqrt{\mathrm{var_A}(f)\mathrm{var_F}(f)}
$
is called \texttt{the  Breitenberger (periodic) UC}. The corresponding  
uncertainty principle \cite{B, PQ} says that if $f \in L_{2,2\pi}$,   $f(x)\neq C 
\mathrm{e}^{ \mathrm{i} k x},$ $C \in \mathbb{R},$ $k \in \mathbb{Z}$, then  $UC_B(f) > 1/2$ and there is no function such that  $UC_B(f) = 1/2$


\section{Results}

\subsection{Nonstationary wavelets generated by periodic wavelets}

\begin{lem}
\label{def_mask}
Let $\nu^j_k\in\mathbb{R},$ $\nu^j_k \geq 0,$ be a two-parametric sequence such that 
$\nu^j_k=\nu^j_{k+2^j},$ $|\nu^j_k|^2 + |\nu^j_{k+2^{j-1}}|^2=1,$ 
$\nu^j_k = \nu^j_{-k}.$ By definition, put 
$\theta^j_k:=\arccos \nu^j_k.$
Let a function $z^K_j$ defined on the interval $[0,\,1/4]$ be a uniform spline of order $K+1$ of minimal defect such that $z^K_j(k 2^{-j}) = \theta^j_k,$ $k=0,\dots,2^{j-2},$ $(z^K_j)^{(l)}(0)=0,$ $l=1,\dots,K-1.$
Finally, let $m^K_j$ be an even $1$-periodic function defined as
\begin{equation}
\label{mask}
m^K_j(\xi)=
\left\{
\begin{array}{ll}
\cos(z^K_j(\xi)), & \xi \in [0,\,1/4], \\
\sin(z^K_j(1/2-\xi)), & \xi \in (1/4,\,1/2].
\end{array}
\right.
\end{equation}

Then $|m^K_j(\xi)|^2+|m^K_j(\xi+1/2)|^2=1,$  $m^K_j \in C^{K-1}(\mathbb{T}),$ 
$(m^K_j)^{(l)}(1/2) = 0,$ $l=1,\dots,K-1.$  
\end{lem}
\textbf{Proof.}

The equality 
$|m^K_j(\xi)|^2+|m^K_j(\xi+1/2)|^2=1$
follows immediately from 
the detailed definition of $m^K_j$ 
$$
m^K_j(\xi)=
\left\{
\begin{array}{ll}
\sin(z^K_j(1/2+\xi-k)), & \xi \in [-1/2+k,\, -1/4+k), \\
\cos(z^K_j(-\xi+k)), & \xi \in [-1/4+k,\, k), \\
\cos(z^K_j(\xi-k)), & \xi \in [k,\, 1/4+k), \\
\sin(z^K_j(1/2-\xi+k)), & \xi \in [1/4+k,\, 1/2+k), k \in \z.
\end{array}
\right.
$$
For instance, for $\xi \in [-1/2+k,\, -1/4+k)$ we get
$$
|m^K_j(\xi)|^2+|m^K_j(\xi+1/2)|^2= \sin^2(z^K_j(1/2+\xi-k)) + \cos^2(z^K_j(1/2+\xi-k)) =1.
$$
Since $z^K_j$ is a spline of order $K+1$ of minimal defect, we obtain $z^K_j \in C^{K-1}(\mathbb{T}).$ Therefore, $m^K_j$ is $K-1$ times continuously differentiable at all points except $n/4,$ where $n\in\z.$ 
The smoothness of  $m^K_j$ at the origin follows from the fact that $\cos$ is smooth and even. Since $m^K_j$ is periodic and even, it remains to check  smoothness at points 
$\xi = 1/4,$ $\xi=1/2.$
At the point $\xi=1/4$, we get for 
$n = 0, \dots, K-1$
$$
\left.(m^K_j)^{(n)}(\xi)\right|_{\xi =1/4-0} =
\left.
\left(\cos \left( z^K_j(\xi)  \right)\right)^{(n)}
 \right|_{\xi =1/4-0} 
$$
$$
\left.(m^K_j)^{(n)}(\xi)\right|_{\xi =1/4+0} =
\left.
\left(\cos \left(\pi/2 - z^K_j(1/2 - \xi)  \right)\right)^{(n)}
 \right|_{\xi =1/4+0}.
$$
Thus, using $z^K_j(1/4)=\theta^j_{2^{j-2}}=\pi/4$, we get $(m^K_j)^{(n)}(1/4-0) = (m^K_j)^{(n)}(1/4+0)$ for $n = 0, \dots, K-1.$
At the point $\xi=1/2$ we have  
$$
\left.(m^K_j)^{(n)}(\xi)\right|_{\xi =1/2-0} = 
\left.(\sin(z^K_j(1/2-\xi)))^{(n)} \right|_{\xi =1/2-0}
$$
$$
\left.(m^K_j)^{(n)}(\xi)\right|_{\xi =1/2+0} = 
\left.(\sin(z^K_j(-1/2+\xi)))^{(n)} \right|_{\xi =1/2+0}
$$ 
So, by the condition $(z^K_j)^{(l)}(0)=0,$ $l=1,\dots,K-1$ we immediately obtain
$$
(m^K_j)^{(n)}(1/2-0) =
(m^K_j)^{(n)}(1/2+0) = 0.
$$
This concludes the proof of Lemma \ref{def_mask}. \hfill $\Diamond$

\begin{lem}
\label{prodL2}
If $a_j\in L_2(\mathbb{T}),$ $a_j(0)=1$, and $\sum_{j}\|a''_j\|_2/2^j<\infty$, then $\widehat{\varphi_j}(\xi)=\prod_{r=1}^{\infty} a_{j+r}(\xi/2^r)$ uniformly and absolutely converges on any $[a,\,b]\subset \mathbb{R}.$ If additionally $|a_j(\xi)|^2+|a_j(\xi+1/2)|^2=1,$ 
then  $\widehat{\varphi_j} \in L_2(\mathbb{R})$ and $\|\widehat{\varphi_j}\|_2\leq 1.$
\end{lem}
 
\textbf{Proof.} 
This is a slight modification of the corresponding stationary result 
(see \cite[Proposition 2.4.1]{NPS}). Suppose $a_j(\xi) = \sum_{k\in \z}c_{j,k} e^{2 \pi i k \xi}$ and $a''_j \in L_2(\mathbb{T})$. Using $a_j(0)=1$, we get
$$
|a_j(\xi) - 1|= \left|\sum_{k\in \z} c_{j,k} \left( e^{2 \pi i k \xi} - 1\right)\right| 
\leq 2\pi \sum_{k\in \z} |c_{j,k}| |k| |\xi|
\leq  C  |\xi|  \left(\sum_{k\in \z} | k^2 c_{j,k}|^2 \right)^{1/2} = 
C_1  |\xi| \|a''_j\|_2.
$$
Therefore,
$$
\sum_{r=1}^{\infty}|a_{j+r}(\xi/2^r) - 1| \leq C_1 |\xi| \sum_{r=1}^{\infty} 
\frac{\|a''_{j+r}\|_2}{2^r} =  C_1 2^j  |\xi| \sum_{n=1}^{\infty} 
\frac{\|a''_{n}\|_2}{2^n}.
$$
Hence, the infinite product 
$\prod_{r=1}^{\infty} a_{j+r}(\xi/2^r)$ uniformly with respect to $\xi$ and absolutely converges on any $[a,\,b]\subset \mathbb{R}.$ The proof of the facts $\widehat{\varphi_j} \in L_2(\mathbb{R})$ and $\|\widehat{\varphi_j}\|_2\leq 1$ can be rewritten from stationary case \cite[Lemma 4.1.3]{NPS}. Lemma \ref{prodL2} is proved. \hfill $\Diamond$

\begin{lem}
\label{my_prodL2}
If $m^K_j$ is defined by (\ref{mask}), $\nu^j_0=1,$ and 
$$
\sum_{j}\frac{1}{2^j}\left(\int_{0}^{1/4}((z^K_j)'(\xi))^4+((z^K_j)''(\xi))^2\,d\xi\right)^{1/2}<\infty,
$$
then 
$\widehat{\varphi_j}(\xi)=\prod_{r=1}^{\infty} m^K_{j+r}(\xi/2^r)$ uniformly and absolutely converges on any $[a,\,b]\subset \mathbb{R},$ 
 $\widehat{\varphi_j} \in L_2(\mathbb{R})$ and $\|\widehat{\varphi_j}\|_2\leq 1.$
\end{lem}

\textbf{Proof.}
Assumptions of Lemma \ref{my_prodL2} are specifications of conditions of Lemma \ref{prodL2} for the masks $m^K_j$.
Indeed, 
$$
\|(m^K_j)''\|^2_2 = \int_0^1 |(m^K_j)''(\xi)|^2 \, d\xi = 
2 \left(\int_0^{1/4} (\cos''(z^K_j(\xi)))^2 \, d\xi 
+ \int_{1/4}^{1/2} (\sin''(z^K_j(1/2 - \xi)))^2 \, d\xi \right)
$$
$$
\ \ \ \ \ \ \ \ \ \ \ \ \ \ \ \ \ \ \ \ \ \ \ \ \ \ \ \ \ \ \ \ \ \ \ \ \  
=
2\left(\int_0^{1/4} ((z^K_j)'(\xi))^4 + ((z^K_j)''(\xi))^2 \, d\xi \right) 
\ \ \ \ \ \ \ \ \ \ \ \ \ \ \ \ \ \ \ \ \ \ \ \ \ \ \ \ \ \ \ \ \ \ \ \ \ \ 
\Diamond
$$


\subsection{Adjustment of localization}

The following theorem describes a connection between $UC_H$ and $UC_B$ in nonstationary case.
For stationary setup, this theorem is proved in \cite{prqurase03}. 

\begin{teo}
\label{nstPQRS}
Let  $\Psi^P =\{\varphi^{P}_{0}, \, \psi^{P}_{j,k}  \}_{j\in \mathbb{Z}_+,\, k=0,\dots,2^j-1}$
and
$\Psi^N = \{\varphi^{N}_{0,k}, \, \psi^{N}_{j,k} \}_{j\in \mathbb{Z}_+,\, k \in \mathbb{Z}}$ be  periodic and nonstationary Parseval wavelet frames, and 
 
$$
\psi^P_j(x):=\sum_{n\in\mathbb{Z}} \psi^N_j (x-n).
$$ 

If there exist  functions $f,\, f_1\in L_2(\mathbb{R})$ such that
$|2^{-j/2}\psi^{N}_j(2^{-j} x)|\leq f(x),$ $|(2^{-j/2}\psi^{N}_j(2^{-j} x))'(x)|\leq f_1(x),$ and $f\in AC_{loc}(\mathbb{R})$, $f(x)=O(|x|^{-3/2-\varepsilon}),$ $f_1(x)=O(|x|^{-1-\varepsilon})$ as $x\to \infty,$ $\varepsilon > 0$, then 

$$
\lim_{j\to \infty} UC_B(\psi^P_j) = \lim_{j\to \infty}UC_H(\psi^N_j).
$$  
\end{teo}

We omit the proof of the theorem since it can be straightforwardly checked that all the steps of the proof of Theorem 3 \cite{prqurase03} holds true for nonstationary case under the assumptions of Theorem  \ref{nstPQRS}.

It is not convenient to apply  Theorem  \ref{nstPQRS}. Indeed, our stating point is a periodic wavelet system, however the main condition (existence of majorants $f$, $f_1$) concerns the resulting nonstationary system. The next theorem is free of this drawback and provides sufficient conditions for an adjustment of localization in terms of initial periodic masks.

\begin{teo}
\label{adjust_loc} 
Let 
$\Psi^P =\{\varphi^{P}_{0}, \, \psi^{P}_{j,k}  \}_{j\in \mathbb{Z}_+,\, k=0,\dots,2^j-1}$ 
be a periodic Parseval wavelet frame, $(\mu^j_k)_k=(2^{1/2} \nu^j_k)_k$ be its  scaling masks. Let $(\nu^j_k)_k$ satisfy the conditions of Lemma \ref{def_mask}. 
 By definition, put $\theta^j_k:=\arccos \nu^j_k,$ 
$\overline{\nu}^j_k:=\max\{\nu^j_k,\,\nu^j_{k+1}\}$,
$\underline{\nu}^j_k:=\min\{\nu^j_k,\,\nu^j_{k+1}\}$. 
If
\begin{enumerate}
 \item
the series $\sum_{k\in\z}|k b^j_k|$ uniformly converges and uniformly bounded with respect to 
$j$ , where
	$b^j_k:=\prod_{r=1}^{\infty} \overline{\nu}^{j+r}_k$, and
\item
 $
|\theta^j_{k+1}-\theta^j_{k}|\leq C 2^{-j},
	$
where $C$ is an absolute constant,
\end{enumerate}

\noindent
 then the system $\Psi^N = \{\varphi^{N}_{0,k}, \, \psi^{N}_{j,k} \}_{j\in \mathbb{Z}_+,\, k \in \mathbb{Z}}$ with scaling masks $m^1_j(\xi)$ defined in Lemma \ref{def_mask} as $K=1$  forms a nonstationary Parseval wavelet frame,
 and 
$
\lim_{j\to \infty} UC_B(\psi^P_j) = \lim_{j\to \infty}UC_H(\psi^N_j).
$  
\end{teo}  

\begin{rem}
Convergence of the infinite products  
$a^j_k:=\prod_{r=1}^{\infty} \underline{\nu}^{j+r}_k$ 
and
$b^j_k:=\prod_{r=1}^{\infty} \overline{\nu}^{j+r}_k$  follows from the  conditions of Theorem \ref{adjust_loc}.
\end{rem}
 
\textbf{Proof.}
We check that under assumptions 1. and 2. all the conditions of Theorem \ref{nstPQRS} are fulfilled. 

First of all, the infinite products   
$a^j_k=\prod_{r=1}^{\infty} \underline{\nu}^{j+r}_k$ 
and
$b^j_k=\prod_{r=1}^{\infty} \overline{\nu}^{j+r}_k$ 
converges for any $j\in\n$ and $k\in\z.$ Infinite products are considered convergent also in the case when they are equal to zero. 
Indeed, it follows from (\ref{con2}) and definition of $\overline{\nu}^{j}_k$ that  $1 \geq \overline{\nu}^{j}_k \geq \nu^j_k \geq 0$.
Therefore, $\bigl|\log \overline{\nu}^{j}_k \bigr| \leq \bigl|\log \nu^j_k \bigr|$, and we know that  the series $\sum_{r=1}^{\infty} \log \nu^j_k$ is absolutely  convergent or equal to $-\infty$. So, the series  
 $\sum_{r=1}^{\infty} \log \overline{\nu}^{j}_k$ is absolutely  convergent or equal to $-\infty$, hence, the product $\prod_{r=1}^{\infty} \overline{\nu}^{j+r}_k$ is  absolutely  convergent. Concerning to the infinite product 
$a^j_k=\prod_{r=1}^{\infty} \underline{\nu}^{j+r}_k$, we note 
that since the product 
$\prod_{r=1}^{\infty} \nu^{j+r}_k$ is  absolutely  convergent, the series 
$\sum_{r=1}^{\infty} (\nu^{j+r}_k - 1)$ has the same property. Therefore, 
using 2. we get 
$$
\sum_{r=1}^{\infty} |\underline{\nu}^{j+r}_k - 1|=
\sum_{r=1}^{\infty} (1-\underline{\nu}^{j+r}_k)=
\sum_{r=1}^{\infty} (1 - \nu^{j+r}_k) +
\sum_{r=1}^{\infty} (\nu^{j+r}_k - \underline{\nu}^{j+r}_k) 
\leq
\sum_{r=1}^{\infty} (1 - \nu^{j+r}_k) + 
C \sum_{r=1}^{\infty} \frac{1}{2^{j+r}},
$$
so the series 
$\sum_{r=1}^{\infty} (\underline{\nu}^{j+r}_k - 1)$ and, hence, the product  
$\prod_{r=1}^{\infty} \underline{\nu}^{j+r}_k$
are  absolutely  convergent.

The next step is to check that the infinite product 
$\prod_{r=1}^{\infty} m^1_{j+r}(\xi/2^r) = \widehat{\varphi_j}(\xi)$ 
is uniformly convergent on any interval $[a,\,b]\in\r$ and 
$\varphi_j \in L_2(\r),$ where the auxiliary mask $m^1_{j}$ is defined in Lemma \ref{def_mask}.   Lemma \ref{my_prodL2} can not help us here since  it works as $K\geq 2$, and we consider the case $K=1.$ We check aforementioned properties directly. It follows from an elementary property of $\cos$ that the mask $m^1_{j}$ can be rewritten as 
$ m^1_j(\xi) = \cos (z^1_j(\xi)).$ The function $z^1_j$ is piecewise linear, namely, 
$
z^1_j(\xi) = 2^j(\theta^j_{k+1}-\theta^j_{k})(\xi - k 2^{-j}) + \theta^j_{k},
$
$\xi \in [k 2^{-j},\, (k+1) 2^{-j}].$
 Therefore, $z^1_j(\xi)$ 
lies between $\theta^j_{k}$ and $\theta^j_{k+1}$ for 
$\xi \in [k 2^{-j},\, (k+1) 2^{-j}].$ Hence, 
$\underline{\nu}^{j}_k\leq  m^1_j(\xi) \leq \overline{\nu}^{j}_k$
and
\begin{equation}
\label{est_phi}
a^j_k\leq \widehat{\varphi_j}(\xi):= \prod_{r=1}^{\infty}  m^1_{j+r}(\xi/2^r) \leq b^j_k
\mbox{ for } \xi \in [k 2^{-j},\, (k+1) 2^{-j}].
\end{equation}
So, the product 
$
\prod_{r=1}^{\infty}  m^1_{j+r}(\xi/2^r) 
$     
is everywhere finite. Using $0 \leq m^1_j(\xi)\leq 1$, we get
 $$
\left| \prod_{r=1}^{n}  m^1_{j+r}(\xi/2^r) - \widehat{\varphi_j}(\xi) \right| \leq
1- \prod_{r=n+1}^{\infty}  m^1_{j+r}(\xi/2^r) \leq 
1-\prod_{r=n+1}^{\infty} \underline{\nu}^{j+r}_k \to 0 
\mbox{ as } n \to \infty.
$$
Since any interval $[a,\,b]$ can be covered by a finite number of the intervals $[k 2^{-j},\,(k+1) 2^{-j}]$, the convergence is uniform on $[a,\,b]$. Using the inequality $0\leq\widehat{\varphi_j}(\xi) \leq b^j_k$ for $\xi \in  [k 2^{-j},\, (k+1) 2^{-j}]$ and the following corollary from condition 1. 
$\sum_{k\in \z} (b^j_k)^2 < \infty $ we immediately get 
$\widehat{\varphi_j} \in L_2(\r),$ thus $\varphi_j \in L_2(\r).$ 

Now we claim that the system $\Psi^N = \{\varphi^{N}_{0,k}, \, \psi^{N}_{j,k} \}_{j\in \mathbb{Z}_+,\, k \in \mathbb{Z}}$ with the scaling masks 
$m^N_{0,j}(\xi) = 2^{1/2} m^1_j(\xi),$ where the auxiliary masks
$m^1_j(\xi)$ are defined in Lemma \ref{def_mask} as $K=1$,  forms a nonstationary Parseval wavelet frame. It follows from Theorem \ref{nUEP} and (\ref{aux_n}). Indeed,  
$\varphi_j \in L_2(\r)$ is already checked, conditions 
(\ref{ncon2}) and (\ref{ncon3}) are provided by Lemma \ref{def_mask} and a definition of $\widehat{\varphi_j}$ by means of  the infinite product. Condition (\ref{ncon1}) 
is equivalent to $\widehat{\varphi_j}(\xi 2^{-j}) \to 1$ as $j \to \infty$ for a fixed $\xi \in \r.$ Taking into account convergence of the infinite products 
$a^j_k=\prod_{r=1}^{\infty} \underline{\nu}^{j+r}_k$
 and 
$b^j_k=\prod_{r=1}^{\infty} \overline{\nu}^{j+r}_k$,
 we obtain $a^j_k \to 1$ and $b^j_k \to 1$ as $j \to \infty.$ Then (\ref{ncon1}) follows from  inequality (\ref{est_phi}).

Switching to property of adjustment of localization, first, we show that for $\xi \neq k 2^{-j}$ 
\begin{equation}
\label{der1}
\widehat{\varphi_j}'(\xi)=\sum_{q=1}^{\infty} 
g_{j,q}(\xi) \left( m^1_{j+q}(\xi/2^q) \right)',
\end{equation}
where $g_{j,q}(\xi):=\prod_{r=1, r\neq q}^{\infty}  m^1_{j+r}(\xi/2^r)$
and
\begin{equation}
\label{der2}
\widehat{\varphi_j}''(\xi)=\sum_{t=1}^{\infty} 
g_{j,t}'(\xi) \left( m^1_{j+t}(\xi/2^t) \right)'+
g_{j,t}(\xi) \left( m^1_{j+t}(\xi/2^t) \right)''
\end{equation}
To prove (\ref{der1}),  it is sufficient to check that
$$
\sum_{q=1}^{n} 
\prod_{r=1, r\neq q}^{n}  m^1_{j+r}(\xi/2^r) \left( m^1_{j+q}(\xi/2^q) \right)' - 
\sum_{q=1}^{\infty} 
g_{j,q}(\xi) \left( m^1_{j+q}(\xi/2^q) \right)' \to 0
$$
as $n \to \infty$ uniformly on any interval $[a,\,b].$ The last difference
is equal to
$$
\left[\sum_{q=1}^{n} 
\prod_{r=1, r\neq q}^{n}  m^1_{j+r}(\xi/2^r) \left( m^1_{j+q}(\xi/2^q) \right)'\right]\left(1-\prod_{r=n+1}^{\infty}  m^1_{j+r}(\xi/2^r)\right)
$$
$$ 
- \sum_{q=n+1}^{\infty} 
g_{j,q}(\xi) \left( m^1_{j+q}(\xi/2^q) \right)'=:
s_1(\xi)-s_2(\xi).
$$
To estimate $s_1(\xi)$, we notice that 
$$
0 \leq \prod_{r=1, r\neq q}^{n}  m^1_{j+r}(\xi/2^r)\leq 1,
$$
\begin{equation}
\label{der_mask}
\left( m^1_{j+q}(\xi/2^q) \right)' = - 2^j (\theta^{j+q}_{k+1}-\theta^{j+q}_{k}) \sin(z^1_{j+q}(\xi/2^q))
\end{equation}
$$
0\leq 1-\prod_{r=n+1}^{\infty}  m^1_{j+r}(\xi/2^r)\leq 
1-\prod_{r=n+1}^{\infty} \underline{\nu}^{j+r}_k.
$$
Combining these estimates and condition 2., we get
$$
|s_1(\xi)|\leq \sum_{q=1}^n 2^j \bigl|\theta^{j+q}_{k+1}-\theta^{j+q}_{k}\bigr|
\left(1-\prod_{r=n+1}^{\infty} \underline{\nu}^{j+r}_k\right) 
\leq C \sum_{q=1}^n \frac{2^j}{2^{j+q}}   \left(1-\prod_{r=n+1}^{\infty} \underline{\nu}^{j+r}_k\right)
\leq C   \left(1-\prod_{r=n+1}^{\infty} \underline{\nu}^{j+r}_k\right) 
$$ 
The last expression tends to zero as $n \to \infty.$

To estimate $s_2(\xi)$, we use consequently  $g_{j,q} \leq 1$, (\ref{der_mask}), and condition 2., so, we have
$$
|s_2(\xi)| \leq \sum_{q=n+1}^{\infty} 2^j \bigl|\theta^{j+q}_{k+1}-\theta^{j+q}_{k}\bigr|\leq C\sum_{q=n+1}^{\infty} \frac{2^j}{2^{j+q}} = \frac{C}{2^n}.
$$
Thus (\ref{der1}) is proved. Switching to proof of (\ref{der2}) we notice that 
$$
g'_{j,t}(\xi)=\sum_{q=1}^{\infty}\prod_{r=1, r\neq q, r\neq t}^{n}  m^1_{j+r}(\xi/2^r) \left( m^1_{j+q}(\xi/2^q) \right)'
$$ 
(It can be checked analogously to (\ref{der1}).) So, 
using $0 \leq m^1_j \leq 1$, (\ref{der_mask}), and condition 2., we have
$$
|g'_{j,t}(\xi)| \leq \sum_{q=1}^{\infty} \left|\left( m^1_{j+q}(\xi/2^q) \right)'\right|  \leq C,
$$  
therefore, again by (\ref{der_mask}), and condition 2.
\begin{equation}
\label{est1}
\sum_{t=n+1}^{\infty} 
\left|g_{j,t}'(\xi) \left( m^1_{j+t}(\xi/2^t) \right)'\right| \leq
C \sum_{t=n+1}^{\infty} 
\left| \left( m^1_{j+t}(\xi/2^t) \right)'\right| \leq 
\frac{C^2}{2^n}.
\end{equation}
Analogously,
$$
\sum_{t=n+1}^{\infty} 
\left|g_{j,t}(\xi) \left( m^1_{j+t}(\xi/2^t) \right)''\right| \leq
\sum_{t=n+1}^{\infty} 
\left| \left( m^1_{j+t}(\xi/2^t) \right)''\right|\leq 
\frac{C^2}{2^{2n}}
$$
Combining the last estimate and (\ref{est1}) we get (\ref{der2}).
Uniform convergence of (\ref{der1}) and (\ref{der2}) implies continuity of $\widehat{\varphi_j}'$ and $\widehat{\varphi_j}''$ on the intervals 
$(k 2^{-j},\, (k+1) 2^{-j}).$

Now we are ready to provide majorants $f,\, f_1 \in L_2(\r)$ such that 
$|\psi_j(x)|\leq f(x),$ $|\psi'_j(x)|\leq f_1(x)$,  $f(x) = O(|x|^{-2})$,
$f_1(x) = O(|x|^{-2})$ as $|x| \to \infty,$
where $\psi_j(x)= 2^{-j/2}\psi^N_j(2^{-j}x)$ are auxiliary wavelet functions.
First, we get majorants for scaling functions $\varphi_j.$ 
Using condition 1. and the inequality $|\widehat{\varphi_j}(\xi)|\leq b^j_k$ for $\xi \in [k 2^{-j},\,(k+1)2^{-j}]$  we get $\widehat{\varphi_j} \in L_1(\r),$  
$\xi \widehat{\varphi_j}(\xi) \in L_1(\r).$
 The product 
$\prod_{r=1}^{\infty}  m^1_{j+r}(\xi/2^r)$ is uniformly convergent on any 
$[a,\,b],$ therefore, the function $\widehat{\varphi_j}$ is continuous on $\mathbb{R}$. Thus 
$$
\varphi_j(x)=\int_{\r}\widehat{\varphi_j}(\xi)e^{2\pi i \xi x}\, d\xi 
= \left.
\frac{1}{2\pi i x} \widehat{\varphi_j}(\xi)e^{2\pi i \xi x}\right|_{\r}- \left.
\frac{1}{(2\pi i x)^2} \sum_{k \in \z} \widehat{\varphi_j}'(\xi)e^{2\pi i \xi x}\right|^{\xi=\frac{k+1}{2^j}-0}_{\xi=\frac{k}{2^j}+0} 
$$
$$
+
\frac{1}{(2\pi i x)^2}
\int_{\r}\widehat{\varphi_j}''(\xi)e^{2\pi i \xi x}\, d\xi
$$
and
$$
\varphi'_j(x)=\int_{\r} 2 \pi i \xi \widehat{\varphi_j}(\xi)e^{2\pi i \xi x}\, d\xi 
= \left.
\frac{1}{ x}  \xi \widehat{\varphi_j}(\xi)e^{2\pi i \xi x}\right|_{\r}- \left.
\frac{1}{2\pi i x^2} \sum_{k \in \z} \left(\xi \widehat{\varphi_j}(\xi) \right)'e^{2\pi i \xi x}\right|^{\xi=\frac{k+1}{2^j}-0}_{\xi=\frac{k}{2^j}+0} 
$$
$$
+
\frac{1}{2\pi i x^2}
\int_{\r}\left(\xi\widehat{\varphi_j}(\xi)\right)''e^{2\pi i \xi x}\, d\xi
$$

 Since $\widehat{\varphi_j} \in L_1(\r)$, $\xi \widehat{\varphi_j}(\xi) \in L_1(\r)$ and 
$\widehat{\varphi_j}'$ is continuous at any point $\xi \neq k 2^{-j}$, it follows that 
 $\widehat{\varphi_j}(\xi) \to 0$ and 
 $\xi \widehat{\varphi_j}(\xi) \to 0$
 as $\xi \to \pm \infty.$  
Therefore,
$$
\left.
\frac{1}{ x}  \xi \widehat{\varphi_j}(\xi)e^{2\pi i \xi x}\right|_{\r}=0,
\ \ \ 
\left.
\frac{1}{ x}  \xi \widehat{\varphi_j}(\xi)e^{2\pi i \xi x}\right|_{\r}=0.
$$
Hence,
\begin{equation}
\label{phi_es}
\varphi_j(x)=- 
\frac{1}{(2\pi i x)^2} \sum_{k \in \z} \left(\widehat{\varphi_j}'\left(\frac{k}{2^j}-0\right)-\widehat{\varphi_j}'\left(\frac{k}{2^j}+0\right)\right)e^{2\pi i \frac{k}{2^j} x}
+
\frac{1}{(2\pi i x)^2}
\int_{\r}\widehat{\varphi_j}''(\xi)e^{2\pi i \xi x}\, d\xi
\end{equation}
\begin{equation}
\label{phi1_es}
\varphi'_j(x)= -\frac{1}{2\pi i x^2} 
\sum_{k \in \z}\frac{k}{2^j} \left(\widehat{\varphi_j}'\left(\frac{k}{2^j}-0\right)-\widehat{\varphi_j}'\left(\frac{k}{2^j}+0\right)\right)e^{2\pi i \frac{k}{2^j} x}
\end{equation}
$$
+
\frac{1}{\pi i x^2}
\int_{\r}\widehat{\varphi_j}'(\xi)e^{2\pi i \xi x}\, d\xi
+\frac{1}{2\pi i x^2}
\int_{\r}\xi\widehat{\varphi_j}''(\xi)e^{2\pi i \xi x}\, d\xi
$$
We denote
$$
A_1:=\left|\sum_{k \in \z} \left(\widehat{\varphi_j}'\left(\frac{k}{2^j}-0\right)-\widehat{\varphi_j}'\left(\frac{k}{2^j}+0\right)\right)e^{2\pi i \frac{k}{2^j} x}\right|, 
$$
$$ 
A_2:= \left|\sum_{k \in \z}\frac{k}{2^j} \left(\widehat{\varphi_j}'\left(\frac{k}{2^j}-0\right)-\widehat{\varphi_j}'\left(\frac{k}{2^j}+0\right)\right)e^{2\pi i \frac{k}{2^j} x}\right|,
$$
$$
A_3:=\left|\int_{\r}\widehat{\varphi_j}'(\xi)e^{2\pi i \xi x}\, d\xi\right|, 
\ \ 
A_4:= \left|\int_{\r}\widehat{\varphi_j}''(\xi)e^{2\pi i \xi x}\, d\xi\right|,
\ \ 
A_5:=\left|\int_{\r}\xi\widehat{\varphi_j}''(\xi)e^{2\pi i \xi x}\, d\xi\right|.
$$
Let us prove that $A_n,$ $n=1,\dots,5$ are uniformly bounded with respect to $j\in\n$ and $x\in\r.$
First, we estimate $\widehat{\varphi_j}'$ and $\widehat{\varphi_j}''$. Let $\xi \in [k 2^{-j},\, (k+1) 2^{-j}].$
Recalling (\ref{der1}), we need to estimate $g_{j,q}(\xi)$ 
and
$(m^1_{j+q}(\xi 2^{-q}))'$.
  Using $0\leq m^1_j(\xi)\leq 1$, definition of $\widehat{\varphi_j},$ and (\ref{est_phi}) we get 
$$
0\leq  g_{j,q}(\xi) \leq \prod_{r=q+1}^{\infty} m^1_{j+r}(\xi / 2^{r}) 
= \widehat{\varphi_{j+q}}(\xi) \leq b^{j+q}_k
$$
By (\ref{der_mask}) and condition 2.
$$
\left|(m^1_{j+q}(\xi 2^{-q}))'\right| \leq 
\left|2^j (\theta^{j+q}_{k+1}-\theta^{j+q}_{k})\right| \leq
\frac{C}{2^q}.
$$
Substituting estimates in (\ref{der1}), we obtain
\begin{equation}
\label{est_phi'}
|\widehat{\varphi_j}'(\xi)|\leq C \sum_{q=1}^{\infty} b^{j+q}_{k} 
\frac{1}{2^q} \mbox{ for }\xi \in [k 2^{-j},\, (k+1) 2^{-j}].
\end{equation}

Now we majorize $\widehat{\varphi_j}''$. We start with (\ref{der2}) and estimate $g'_{j,t}(\xi)$ and $(m^1_{j+t}(\xi 2^{-t}))''$. 
Analogously to (\ref{est_phi'}) we have
$$
|g'_{j,t}(\xi)| \leq 
\sum_{q=1}^{\infty}\prod_{r=1, r\neq q, r\neq t}^{n}  m^1_{j+r}(\xi/2^r) \left|\left( m^1_{j+q}(\xi/2^q) \right)'\right|
\leq  C \sum_{q=1}^{\infty} b^{j+q+t}_{k} 
\frac{1}{2^q}
$$
By (\ref{der_mask}) and condition 2.
$$
\left|(m^1_{j+t}(\xi 2^{-q}))''\right| \leq 
\left|2^j (\theta^{j+t}_{k+1}-\theta^{j+t}_{k})\right|^2 \leq
\frac{C^2}{2^{2t}}.
$$
Collecting all inequalities and substituting them in (\ref{der2}) we obtain
\begin{equation}
\label{est_phi''}
\left|\widehat{\varphi_j}''(\xi)\right| \leq
C^2 \sum_{t=1}^{\infty} \sum_{q=1}^{\infty} b^{j+q+t}_k \frac{1}{2^{q+t}}
+
C^2 \sum_{t=1}^{\infty} b^{j+t}_k \frac{1}{2^{2t}}
\mbox{ for }\xi \in [k 2^{-j},\, (k+1) 2^{-j}].
\end{equation}

Now we return to the proof of boundedness of $A_n$, $n=1,\dots,5.$

Using (\ref{est_phi'}), we get
$$
A_1 \leq
\sum_{k \in \z} \left|\widehat{\varphi_j}'\left(\frac{k}{2^j}-0\right)-\widehat{\varphi_j}'\left(\frac{k}{2^j}+0\right)\right| 
\leq
\sum_{k \in \z} \sum_{t =1}^{\infty} 
b^{j+t}_k 2^j \frac{2C}{2^{j+t}}.
$$
Since 
$
\sum_{k \in \z} b^{j}_k \leq b^j_0 + \sum_{k \in \z}  |k| b^{j}_k 
\leq  1 + \sum_{k \in \z}  |k| b^{j}_k,
$
condition 1. make it possible to change the order of summation, thus
we have 
$$
A_1 \leq 2C  \sum_{t =1}^{\infty} \sum_{k \in \z}
b^{j+t}_k  \frac{1}{2^{t}} \leq 2C^2 \sum_{t =1}^{\infty}\frac{1}{2^{t}} 
\leq 2 C^2. 
$$
Analogously, for $A_2$ we get
$$
A_2 
\leq
\sum_{k \in \z} \sum_{t =1}^{\infty} 
|k| b^{j+t}_k \frac{2C}{2^{j+t}} \leq \frac{2C^2}{2^j}.
$$
Again, using (\ref{est_phi'}) we have
$$
A_3 \leq \int_{\r} \left|\widehat{\varphi_j}'(\xi)\right|\, d\xi
\leq 
\sum_{k \in \z}\int_{k 2^{-j}}^{(k+1) 2^{-j}} \left|\widehat{\varphi_j}'(\xi)\right|\, d\xi 
\leq C
\sum_{k \in \z}\int_{k 2^{-j}}^{(k+1) 2^{-j}} 
\sum_{q=1}^{\infty} b^{j+q}_k \frac{1}{2^q} \, d\xi 
= C \frac{1}{2^j}
\sum_{k \in \z} 
\sum_{q=1}^{\infty} b^{j+q}_k \frac{1}{2^q}
$$ 
$$
= C \frac{1}{2^j}
\sum_{q=1}^{\infty} 
\sum_{k \in \z} 
b^{j+q}_k \frac{1}{2^q} \leq \frac{C^2}{2^j}.
$$
By (\ref{est_phi''}) and condition 1.,
$$
A_4 \leq \int_{\r} \left|\widehat{\varphi_j}''(\xi)\right|\, d\xi
\leq 
\sum_{k \in \z}\int_{k 2^{-j}}^{(k+1) 2^{-j}} \left|\widehat{\varphi_j}''(\xi)\right|\, d\xi 
\leq
C^2 \frac{1}{2^j} \sum_{k \in \z} \left(\sum_{t=1}^{\infty} \sum_{q=1}^{\infty} b^{j+q+t}_k \frac{1}{2^{q+t}}
+ \sum_{t=1}^{\infty} b^{j+t}_k \frac{1}{2^{2t}} \right)
$$
$$
= C^2 \frac{1}{2^j}  \left(\sum_{t=1}^{\infty} \sum_{q=1}^{\infty} \sum_{k \in \z} b^{j+q+t}_k \frac{1}{2^{q+t}}
+ \sum_{t=1}^{\infty} \sum_{k \in \z} b^{j+t}_k \frac{1}{2^{2t}} \right)\leq 
\frac{4 C^3}{3} \frac{1}{2^j}, 
$$
$$
A_5 \leq \int_{\r} \left|\xi \widehat{\varphi_j}''(\xi)\right|\, d\xi
\leq 
\sum_{k \in \z}\int_{k 2^{-j}}^{(k+1) 2^{-j}} \left|\xi\widehat{\varphi_j}''(\xi)\right|\, d\xi 
\leq
 \frac{C^2}{2^j} \sum_{k \in \z} (|k|+1)\left(\sum_{t=1}^{\infty} \sum_{q=1}^{\infty} b^{j+q+t}_k \frac{1}{2^{q+t}}
+ \sum_{t=1}^{\infty} b^{j+t}_k \frac{1}{2^{2t}} \right)
$$
$$
= C^2 \frac{1}{2^j}  \left(\sum_{t=1}^{\infty} \sum_{q=1}^{\infty} \sum_{k \in \z}(|k|+1) b^{j+q+t}_k \frac{1}{2^{q+t}}
+ \sum_{t=1}^{\infty} \sum_{k \in \z} (|k|+1) b^{j+t}_k \frac{1}{2^{2t}} \right)\leq 
\frac{8 C^3}{3} \frac{1}{2^j}. 
$$
In all above estimates changing the order of summation is justified by condition 1. So, all the expressions $A_n,$ $n=1,\dots,5$ are bounded by an absolute constant. Thus, (\ref{phi_es}), (\ref{phi1_es}) yield  
$$
\varphi_j(x) \leq C |x|^{-2}, \ \ \ \varphi'_j(x) \leq  C|x|^{-2}.
$$
On the other hand,
$
\varphi_j(x), \varphi'_j(x)
$ 
are bounded by an absolute constant. Indeed, let $v=0$ or $v=1$. Using (\ref{est_phi}) and condition 1., we obtain 
$$
|\varphi^{(v)}_j(x)| \leq \int_{\r} \left| (2 \pi \xi)^v\widehat{\varphi_j}(\xi)\right| \, d\xi 
\leq
\sum_{k \in \z} \int_{k 2^{-j}}^{(k+1) 2^{-j}}\left|(2 \pi \xi)^v\widehat{\varphi_j}(\xi)\right| \, d\xi \leq 
\sum_{k \in \z} \frac{\bigl(2 \pi (|k|+1)\bigr)^v}{2^j} b^j_k \leq C. 
$$

Thus, for functions $\varphi_j$ and $\varphi_j'$ we provide majorants of the form 
$$
\left\{
\begin{array}{cc}
C, & |x|\leq 1; \\
C/|x|^2, & |x|\geq 1. 
\end{array}
\right.
$$ 

Majorants for auxiliary wavelet functions $\psi_j$ and its derivatives $\psi_j'$ can be obtained analogously. Indeed, one can start with the equalities analogous to the case of scaling functions
$$
\psi_j(x)=\int_{\r}\widehat{\psi_j}(\xi)e^{2\pi i \xi x}\, d\xi 
= \left.
\frac{1}{2\pi i x} \widehat{\psi_j}(\xi)e^{2\pi i \xi x}\right|_{\r}- \left.
\frac{1}{(2\pi i x)^2} \sum_{k \in \z} \widehat{\psi_j}'(\xi)e^{2\pi i \xi x}\right|^{\xi=\frac{k+1}{2^j}-0}_{\xi=\frac{k}{2^j}+0} 
$$
$$
+
\frac{1}{(2\pi i x)^2}
\int_{\r}\widehat{\psi_j}''(\xi)e^{2\pi i \xi x}\, d\xi
$$
and
$$
\psi'_j(x)=\int_{\r} 2 \pi i \xi \widehat{\psi_j}(\xi)e^{2\pi i \xi x}\, d\xi 
= \left.
\frac{1}{ x}  \xi \widehat{\psi_j}(\xi)e^{2\pi i \xi x}\right|_{\r}- \left.
\frac{1}{2\pi i x^2} \sum_{k \in \z} \left(\xi \widehat{\psi_j}(\xi) \right)'e^{2\pi i \xi x}\right|^{\xi=\frac{k+1}{2^j}-0}_{\xi=\frac{k}{2^j}+0} 
$$
$$
+
\frac{1}{2\pi i x^2}
\int_{\r}\left(\xi\widehat{\psi_j}(\xi)\right)''e^{2\pi i \xi x}\, d\xi.
$$
Then using the definition of an auxiliary wavelet sequence 
$
\widehat{\psi_j}(\xi) = e^{\pi i \xi} \overline{m^1_{j+1}(\xi/2+1/2)} \widehat{\varphi_{j+1}}(\xi/2) 
$
and boundedness of 
$m^1_{j+1}(\xi/2+1/2),$ $(m^1_{j+1}(\xi/2+1/2))',$ and $(m^1_{j+1}(\xi/2+1/2))''$ we come to the case of scaling sequences. Theorem \ref{adjust_loc} is proved.\hfill $\Diamond$  

Analyzing the assumptions of Theorem \ref{adjust_loc}, it is easy to see that 
condition 2. means the boundedness of the first divided difference 
$2^j (\theta^j_{k+1} - \theta^j_{k})$ for the data points 
$(k 2^{-j},\, \theta^j_{k}),$ $k\in \z.$  

Unfortunately, Theorem \ref{adjust_loc} is not applicable to the Parseval periodic wavelet frame constructed in \cite{LebPres14}. More precisely, condition 2. is not fulfilled. Indeed, for $k=2^{j-2}$ we get $|\theta^j_{k+1}-\theta^j_{k}| \geq 1/\sqrt{2} - \varepsilon_j,$ where $\varepsilon_j \to 0$ as $j\to \infty.$ Moreover, it follows from the last inequality that for any scaling mask $m^N_j$ of a nonstationary wavelet system that corresponds to the periodic wavelet frame constructed in \cite{LebPres14}  there exists a point $\overline{\xi_j}$ such that
$|(m^N_j)'(\overline{\xi_j})|\geq C 2^j$. 
To construct a periodic wavelet sequence $\psi^P_j$ satisfying assumptions of Theorem \ref{adjust_loc} and the equality $\lim_{j\to \infty} UC_B(\psi^P_j) =3/2$  is a task for future investigation.   

Finally, it is interesting to note that we can always construct a trivial  nonstationary wavelet frame starting with a periodic one. 
Suppose 
$\Psi^P =\{\varphi^{P}_{0}, \, \psi^{P}_{j,k}  \}_{j\in \mathbb{Z}_+,\, k=0,\dots,2^j-1}$ 
is a periodic Parseval wavelet frame, $\mu^j_k=2^{1/2} \nu^j_k$ is  scaling mask of this frame.
Let us define a nonstationary scaling mask as a step function such that 
$$
m^N_j(\xi)=2^{1/2} m_j(\xi)=\mu^j_k, \quad   \xi\in[k/2^k,\,(k+1)/2^j).
$$
Then  
a nonstationary scaling function is the following step function
$$
\widehat{\varphi^N_j}(\xi) 
= 2^{-j/2}\prod_{r=1}^{\infty} m_{j+r}(\xi/2^{j+r})
= \widehat{\varphi^P_j}(k), \ \ \   \xi\in[k/2^k,\,(k+1)/2^j).
$$
Since 
$
\|\widehat{\varphi^N_j}\|^2_2 = \|\varphi^P_j\|^2_2$, the scaling function
$\varphi^N_j$ 
 is in $L_2(\mathbb{R}).$ It generates a nonstationary Parseval wavelet frame.
However, since $\widehat{\varphi^N_j}$ is discontinuous,  $\varphi^N_j$ has poor time localization.

\section*{Acknowledgments}
We great fully acknowledge the funding by the RFBR, grant \#15-01-05796, by Saint Petersburg State University, grant  \#9.38.198.2015, and by Volkswagen Foundation.

\end{document}